\documentclass[12pt]{amsart}
\theoremstyle{definition}
\newtheorem{theorem}{Theorem}[section]
\newtheorem{corollary}[theorem]{Corollary}

\newtheorem{lemma}[theorem]{Lemma}
\newtheorem{proposition}[theorem]{Proposition}

\newtheorem{remark}[theorem]{Remark}
\newtheorem{example}[theorem]{Example}

\newtheorem{definition}[theorem]{Definition}
\newtheorem{question}[theorem]{Question}
\theoremstyle{remark}
\numberwithin{equation}{section}



\begin{document}
\title{The dual of the Bourgain-Delbaen space}
\author{Dale Alspach}
\address{Oklahoma State University \\ Department of Mathematics\\ 401
Mathematical Sciences\\ Stillwater, OK 74078-1058}
\email{alspach@math.okstate.edu}
\subjclass{46B20}
\keywords{$\ell_1$-predual, $\mathcal L_\infty$-space, ordinal index,
Szlenk index, uniform homeomorphism}
\begin{abstract}
It is shown that a $\mathcal L_\infty$-space with separable dual
constructed by Bourgain and Delbaen has small Szlenk index and thus does
not have a quotient isomorphic to $C(\omega^\omega)$. It follows that this
is a $\mathcal L_\infty$-space which is the same size as $c_0$
in the sense of the Szlenk index but does
not contain $c_0.$ This has some consequences in the theory of uniform
homeomorphism of Banach spaces.
\end{abstract}
\maketitle
\bigskip
\section{Introduction}
In 1980 Bourgain and Delbaen \cite{BD}
published a method of constructing $\mathcal
L_\infty$-spaces which produced examples with surprising properties.
At the time one of the most interesting aspects of these spaces was
that they were the first examples of a separable space with
the Radon-Nikodym Property but not
isomorphic to a subspace of a separable dual space. In this paper we are
not concerned with this property of the examples, but instead with the fact
that these $\mathcal L_\infty$-spaces fail to contain $c_0$
and thus cannot be isomorphic to an
isometric $L_1(\mu)$-predual. (See \cite{JZ}.) Such spaces are not well
understood and potentially provide a source of interesting examples.

One of our motivations for
considering these spaces was that in \cite{JLS} it was shown that
a $\mathcal L_\infty$ space with $C(\omega^\omega)$ as a quotient is not
uniformly homeomorphic to $c_0$. Thus a natural question is whether that
means that the only $\mathcal L_\infty$-space which is uniformly
homeomorphic to $c_0$ is $c_0$ itself. One consequence of the results
proved here is to show that there is more work to be done by showing that
there are $\mathcal L_\infty$-spaces other than $c_0$ which fail to
have $C(\omega^\omega)$ as a quotient.

If the parameters in the construction in \cite{BD} are
chosen properly, the dual of the space constructed
is separable and therefore by
\cite{LS} is isomorphic to $\ell_1$. 
Our interest is in the $\text{w}^*$-topology on
$\ell_1$ induced by the example. Because the example does not contain $c_0$
it is clear that this $\text{w}^*$-topology is much different than that induced by
a space such as $C(\alpha)$, $\alpha<\omega_1,$ or by a space of affine
functions. One difficulty is that because the dual is only isomorphic to
$\ell_1$, the standard unit vector basis of $\ell_1$ may not be contained
in the extreme points of the unit ball. This property of isometric
$\ell_1$-preduals is heavily (and often implicitly) used in many analyses of
specific $\mathcal L_\infty$-spaces, e.g., \cite{A3},\cite{A4}.
Thus some replacement for this approach is
necessary. Also the definition of the example is given by constructing
embeddings of finite dimensional $\ell_\infty$-spaces and thus infinite
dimensional information must be extracted from this finite dimensional
presentation.

Our approach is to work with the $\text{w}^*$-closure of the $\ell_1$-basis of the
dual space as an
image of a certain associated
compact space with a convenient structure. The $\text{w}^*$-closure
of the $\ell_1$-basis, $C$, is large enough to contain most of the important
information about the dual, since $D=\overline{\text{co }\pm C}^{\|\cdot\|}$
will contain a multiple of the unit ball. On the other hand we do not have
good information about the extreme points and the $\text{w}^*$-topology of this
set $D$. To
overcome this problem we create this associated compact space and we work
through the Choquet theorem and use special information about $C$ which is
encoded in the associated compact space.

In the next section we will recall the definition of the example as given
in \cite{BD} and we will show that the natural coordinate functionals are a
basis for the dual and are equivalent to the usual unit vector basis of
$\ell_1$. In Section~\ref{sz} we develop an approach to computing the Szlenk
index which allows us to move from information about a subset of the dual
to its signed convex hull. This approach may be useful for estimating the
Szlenk index in other situations and thus we develop the ideas in a fairly
general setting. As part of this we introduce a notion of integration for
ordinal-valued functions of a real variable.
In Section~\ref{szcomp} we estimate the Szlenk
index for each $\epsilon>0.$ In the last section we discuss some possible
extensions of the method of construction given by Bourgain and Delbaen.

We use standard notation and terminology from Banach space theory as may be
found in the books \cite{LTI} and \cite{LTII}. We consider only Banach
spaces over the real numbers although much can be adapted to the complex
case. In Section~\ref{sz} we will
need the Szlenk index, \cite{Szl}, so we recall the definition here. 

\begin{definition}
Let $X$ be a
Banach space and let $A\subset X$ and let $B\subset X^*$. Given
$\epsilon>0$ we define a family of subsets of $B$ indexed by the ordinals
less than or equal to $\omega_1$.

Let  $P_0(\epsilon,A,B)=B.$ If $P_\alpha(\epsilon,A,B)$ has been defined,
let 
\begin{multline}
P_{\alpha+1}(\epsilon,A,B)=\{b\in B:\text{there exist }(a_n)\subset A,
(b_n)\subset P_\alpha(\epsilon,A,B)\\
\text{ such that } \text{w}^*\lim
b_n=b,\lim b_n(a_n)\geq \epsilon, \text{w}\lim a_n=0\}.
\end{multline}
If $\alpha$ is a limit ordinal,
\[P_{\alpha}(\epsilon,A,B)=\cap_{\beta<\alpha} P_{\beta}(\epsilon,A,B).
\]
Let $\eta(\epsilon,A,B)$ be the smallest ordinal $\alpha$ such that
$P_{\alpha}(\epsilon,A,B)=\emptyset.$
\end{definition}
Usually $B$ is a $\text{w}^*$-closed subset of $B_{X^*}$ and $A$ is 
$B_X$, where $B_{X^*}$ and $B_X$ are
the unit balls of $X^*$ and $X$, respectively. If $X^*$ is
separable, then $\eta(\epsilon,A,B)$ is defined and countable. Otherwise
the convention is to define $\eta(\epsilon,A,B)=\omega_1$ if there is no
countable ordinal for which the set $P_\alpha(\epsilon,A,B)$ is empty.
In this paper we will always assume that $A=B_X$ so we will
omit this from the notation and write $P_\alpha(\epsilon,B).$

In the case $A=B_X$ and $X^*$ separable it is often convenient to
use a different definition of the Szlenk
index, which yields a slightly different dependence on $\epsilon$, but in
most applications gives equivalent results. In this case the definition of
$P_{\alpha+1}(\epsilon,A,B)$ is replaced by
\begin{multline}
P_{\alpha+1}(\epsilon,A,B)=\{b\in B:\text{there exists }
(b_n)\subset P_\alpha(\epsilon,A,B)\\
\text{ such that } \text{w}^*\lim
b_n=b,\text{ and for all }n\neq m,\|b_n-b_m\|\geq \epsilon\}.
\end{multline}
We will refer to this second version of the Szlenk index as the modified
Szlenk index.

\section{The Bourgain-Delbaen spaces}
\label{bd}
In this section we describe the construction of $\mathcal L_\infty$-spaces
due to Bourgain and Delbaen. We will depart slightly from their notation
and construction, but this is only a matter of convenience. The approach is
to  build a subspace of $\ell_\infty$ by defining a family of
consistent embeddings of $\ell_\infty^{d_n}$ into $\ell_\infty,$ where
$(d_n)$ is some sequence of integers tending to infinity rapidly. The
sequence $(d_n)$ is defined inductively as are the embeddings.

Fix two positive real numbers $a,b$ and a number $\lambda>1$ such that
$b<a\leq 1$ and $a+2b\lambda<\lambda$.
We define $d_1=1,d_2=2$ and assume that $d_k$ has been defined for $k=1,2,\dots,
n$. We define $d_{n+1}-d_n$ to be the cardinality of the set of tuples
$(\sigma',i,m,\sigma'',j)$ such that $1\leq m <n$, $1\leq i\leq d_m$,
$1\leq j\leq d_n$ and $\sigma'$ and $\sigma''$ are $1$ or $-1.$
By enumerating the set of tuples by the integers $k$, $d_n<k\leq d_{n+1}$, we
can inductively
define a map $\phi$ from $\mathbb N\setminus \{1,2\}$ to the set of such
tuples, $(\sigma',i,m,\sigma'',j).$

For each $k\in \mathbb N$, let $e_k^*$ denote the
$k$-th coordinate functional of $\ell_\infty,$ and $e_k$ the $k$-th
coordinate element, i.e., the element of $\ell_\infty$ which is $0$ at each
coordinate except the $k$-th and $1$ in the $k$-th.
To define the embeddings, let $E_n=[e_k:k\leq d_n]$ for each $n$ and define
for $m<n$
inductively $i_{m,n}:E_m\rightarrow E_n$ as follows.
We define $i_{1,2}(te_1)=te_1=e_1^*(t e_1)e_1$ for all $t$ and suppose that 
$i_{m,n}$ has been defined for all $m<n.$
To define an extension map
from $E_n$ into $E_{n+1}$ for each $k$, $d_n<k\leq
d_{n+1},$ we define a functional $f_{\phi(k)}\in E_n^*$ by
\[f_{\phi(k)}(x)=a \sigma' e_i^*(x)+b \sigma'' e_j^*(x-i_{m,n}\pi_m x),\]
where $\pi_m:\ell_\infty\rightarrow E_m$ is standard projection and
$\phi(k)=(\sigma',i,m,\sigma'',j).$ Then
\[i_{n,n+1}(x)=x+\sum_{k=d_n+1}^{d_{n+1}} f_{\phi(k)}(x)e_k\]
for all $x\in E_n.$ Using this map we can define
$i_{m,n+1}(x)=i_{n,n+1}(i_{m,n}(x))$ for all $m<n$ and $x\in E_m.$
In [BD] it is shown that  $\|i_{m,n}\|\leq \lambda$ for all $m<n,$ and thus
considering $\ell_\infty$ as the dual of $\ell_1,$
the $\text{w}^*$-operator
limit $P_m$ of $(i_{m,n}\pi_m)_{n=m+1}^\infty$ exists for
each $m$. ($P_m(x)$ is just the coordinate-wise limit of $i_{m,n}(x)$ for
each $x$ and each coordinate is eventually constant.)
Notice that we can now replace the definition of $f_{\phi(k)}$
by
\[f_{\phi(k)}(x)=a \sigma' e_i^*(x-P_0 x)+b \sigma'' e_j^*(x-P_m x),\]
where $P_0=0.$
Rewriting this in dual form we have 
\[f_{\phi(k)}(x)=a \sigma' (I-P_0^*)e_i^*(x)+b \sigma'' (I-P_m^*)e_j^*(x).\]

We are interested in the spaces $X_{a,b}=[P_m(E_m):m\in \mathbb N],$ where
$a,b$ are fixed constants as above. It
follows easily that $X_{a,b}$ is a $\mathcal L_\infty$-space and in [BD]
some of the Banach space properties of these spaces are determined. If
$a=1$ the dual of $X_{a,b}$ is non-separable and thus is not of interest to
us here. {\bf Thus we assume that $a<1$ unless otherwise noted.} We will
also suppress the subscripts $a,b$ from now on.

Our first task is to show that the dual of $X$ is isomorphic to
$\ell_1$ in a very
concrete sense. Notice that for each $m$, $P_m$ can be considered either as a
map from $\ell_\infty$ into $X$ or as a map from $X$ into itself. Thus
the range of $P_m^*$ is contained in $[e_k^*:k\leq d_m],$ either in
$\ell_\infty^*$ or by restriction to $X$, as elements of $X^*$.

\begin{proposition}\label{quo}
Let $Q$ be the quotient map from $\ell_\infty^*$ onto $X^*.$ Then $(Q(e_n^*))$
is equivalent to the standard unit vector basis of $\ell_1$ and \newline
$Q[e_n^*:n\in
\mathbb N]=X^*.$
\end{proposition}
\begin{proof} 
Because $\|P_m\|\leq \lambda$ and for $g\in[e_1,e_2,\dots,e_{d_m}]$ and
$k\leq d_m$,
\[
P_m^*Q(e_k^*)(g)=e_k^*(P_m(g))=e_k^*(g),\]
for each $m$, it follows that 
\[\|Q\sum_{k=1}^{d_m} a_k e_k^*\|_{X^*}
\geq \lambda^{-1}\|P_m^*Q(\sum_{k=1}^{d_m} a_k e_k^*)\|_{\ell_\infty^*}=
\lambda^{-1}\|\sum_{k=1}^{d_m} a_k e_k^*\|_{\ell_\infty^*}.\]
This proves the first assertion.

For the second we will show that the $\text{w}^*$-closure of $\{Q(e_n^*):n\in
\mathbb N\}$ is contained in $[Q(e_n^*):n \in \mathbb N].$ It then follows
from the Choquet theorem and Smulian's theorem that $[Q(e_n^*):n\in \mathbb
N]$ is $\text{w}^*$-closed and hence equal to $X^*.$ (See \cite{A4}, Lemma 1.)

Let $x^*$ be a $\text{w}^*$-limit point of $(Qe_k^*)_{k\in M}$, for some infinite
subset $M$ of $\mathbb N$.
We may assume that 
$\lim_{k\in M} Q e_k^*(P_m(e_r))=x^*(P_m(e_r))$ for
each $r\leq d_m$ and each $m$.
Let $\phi(k)=(\sigma_k',i_k,m_k,\sigma''_k,j_k).$  We may also assume, by
passing to a smaller index set if necessary, that
$\sigma'_k=\sigma'$ and $\sigma'_k=\sigma''$ for all $k\in M$. Consider
$(m_k)$. If $\sup m_k=\infty$, then 
$b\sigma''(I-P_{m_k}^*)e_{j_k}^*(x)=0$  for all $x\in P_m(E_m)$ for $m\leq m_k$
and thus any $\text{w}^*$-limit point
of $(e_k^*)_{k\in M}$ is a $\text{w}^*$-limit point
of $(a\sigma'(I-P_0^*)e_{i_k}^*).$ 
If $\sup m_k=m<\infty$, then $i_k\leq d_m$ and
 $(a\sigma'(I-P_0^*)e_{i_k}^*)$ has a constant subsequence. Thus any
$\text{w}^*$-limit point of $(e_k^*)_{k\in M}$ is of the form
$a\sigma'(I-P_0^*)e_{i_k}^*+y^*$ where $y^*$ is a $\text{w}^*$-limit point of
$(b\sigma''(I-P_{m_k}^*)e_{j_k}^*).$ Notice that in both cases we have
replaced looking for a
$\text{w}^*$-limit of $(e_k^*)=((I-P_0^*)e_k^*)$ by looking
for a $\text{w}^*$-limit of $(c(I-P_{m_k}^*)e_{r_k}^*)$
where $|c|=a\text{ or }b.$
Therefore we can find a convergent (absolutely summable)
series of terms of the form $c_j(I-P_{m_j}^*)e_{r_k}^*$, $|c_j|\leq a^{j-1}$,
 with limit $x^*$. Actually $c_j=\pm a^sb^{j-s}$ for some $s$, $0\leq s\leq
j,$ and $c_{j+1}=\pm a
c_j$ or $c_{j+1}=\pm b c_j.$ Because $(I-P_m^*)(e_k^*)\in [e_j^*:j\in
\mathbb N]$, for all $m,k$, it follows that $x^*\in [e_j^*:j \in \mathbb N].$
\end{proof}

\begin{remark} In \cite{GKL,GKL1} it is shown that a Banach space which is
uniformly homeomorphic to $c_0$ must have Szlenk index which behaves as
the Szlenk index of $c_0$. It may be possible to use the representation of
the $\text{w}^*$-closure of
the $\ell_1$-basis contained in the previous proof to
get a lower estimate on the Szlenk index and thereby
show that the Bourgain-Delbaen space is not uniformly homeomorphic
to $c_0$.
\end{remark}

\section{Estimating Ordinal Indices}
\label{sz}
We begin by considering an abstract system of derived sets of a metric
space. Eventually we will consider the specific cases where this is the
usual topological derived sets or the Szlenk sets.

\begin{definition} Let $K$ be a closed
subset of a topological space $(X,\tau)$ and let $d(\cdot,\cdot)$ be a
metric on $X$ (which may not be compatible with the topology $\tau$).  A
$\delta${\it -system of derived sets} is
a family, $(K^{(\alpha)})_{\alpha<\omega_1}$, of closed subsets of $K$ such
that
\begin{enumerate}
\item{} there exists some ordinal $\beta_0<\omega_1$ such that
$K^{(\alpha)}=\emptyset$ if $\alpha>\beta_0,$
\item{} if $\alpha<\beta$, then $K^{(\alpha)}\supset K^{(\beta)}$,
\item{} if $\beta$ is a limit ordinal, $\cap_{\alpha<\beta}
K^{(\alpha)}=K^{(\beta)}.$
\item{} 
if $x_n \in K^{(\alpha)}$ for all $n \in \mathbb N$ and  $d(x_n,x_m)\geq
\delta$ for all $n\neq m, n,m\in \mathbb N$ and $\tau-\lim x_n=x$, then $x\in
K^{(\alpha+1)}.$
\end{enumerate}
\end{definition}

For each $\alpha<\omega_1$ let $K^{d(\alpha)}
=K^{(\alpha)}\setminus K^{(\alpha+1)}.$

We are interested in determining how a Szlenk-like index of a set of
finite positive
measures on $K$ as elements of $C(K)^*$ behaves with
respect to this derivation on $K$. To measure this we introduce for each
$\epsilon>0$ and  finite measure  $\mu$ on $K$ the $\epsilon$-distribution
function of $\mu$, 
$f_{\epsilon,\mu}$, from $(0,\infty)$ into $[0,\omega_1)$ 
but with support in
$(0,\epsilon]$. 

To understand the approach consider the following problem. Suppose that $g$
is a nice function on $(0,\infty)$ with values in the countable ordinals. Is
there a sensible notion of area under the graph of $g$?

Because it is not at all clear how to multiply real numbers and ordinals,
let's take a discrete approach. Fix $\epsilon>0.$ For an indicator function
$\gamma 1_{(0,n\epsilon)}$ where $n\in \mathbb N$, we want the
$\epsilon$-area to be $\gamma\cdot n$.
Given an ordinal valued function $g$ on
$(0,\infty)$ the $\epsilon$-area under $g$  should be
the supremum of the ordinal sums
$\gamma_1+\dots+\gamma_k$ of $\epsilon$-areas of disjoint
$\epsilon$-``rectangles'' of width
$\epsilon$ and height $\gamma_i$, $i=1,2,\dots,k,$ which  fit under the
graph of $g$. By a ``rectangle'' we mean a set of the form $A\times B$
where $A$ is Lebesgue measurable and 
$B$ is an interval. There is another difficulty in this in that the
non-commutativity of the addition makes this sensitive to the order in
which the rectangles are taken. To control this difficulty we need the order
of the addition of the rectangles to reflect the values of the function
$g$. To deal with this we use a geometric approach. We think of the ordinal sum
$\gamma_1+\dots+\gamma_k$ as the value of a new function $g'$ on 
$(0,\epsilon]$ 
with $\epsilon$-area under $g'$ approximating  the $\epsilon$-area under $g$.
To be an admissible approximation we require that for each $x$ the segments
in the rectangles above $x$ be in an order which respects the order of the
corresponding segments  under the graph of $g$. More precisely, there is an
injective
function $\psi$ from $\{(x,y):0<x\leq\epsilon,0\leq y \leq g'(x)\}$ into
$\{(x,y):0<x,0\leq y \leq g(x)\}$ such that if for some $x$,
$\psi(x,y_1)=(s,t_1)$ and
$\psi(x,y_2)=(s,t_2)$, then $t_1<t_2$ implies $y_1<y_2.$ Thus the region
under $g'$ is the image under an order preserving (in the second
coordinate only) rearrangement of a
portion of the region under $g$.

At first it may seem that we have drifted far from the original problem.
The connection to our problem is that intuitively the $\epsilon$-Szlenk
index does something similar to computing the
$\epsilon$-area under the distribution of
a measure. Before we introduce precise formulations, consider the measure 
\[\mu = \frac{3}{4}\delta_\omega +
\frac{1}{4}\delta_{\omega^\omega}\]
in the dual of $C(\omega^\omega)$  and its position in the Szlenk sets of
the ball of $C(\omega^\omega)^*$.
Notice that if $1/2<\epsilon\leq 3/4$, $\mu$
is in $P_1(\epsilon)$ but no higher Szlenk set. If $1/4<\epsilon\leq 1/2$,
$\mu$ is in $P_2(\epsilon)$, and if $\epsilon \leq 1/4$, $\mu$ is in
$P_{\omega+3}(\epsilon).$ 
Now consider the distribution function 
\[g(t)=\omega\text{\large 1}_{(0,1/4]}+\text{\large 1}_{(1/4,1]}\] 
and notice that the $\epsilon$-area we have loosely defined above
is the same as the $\epsilon$-Szlenk index of $\mu$, i.e., the $3/4$-area 
is $1$, the $1/2$-area is $2$ and the $1/4$-area is $\omega+3.$

Now we will begin making these ideas precise.
The definition of the $\epsilon$-distribution 
function is via an inductive procedure. We will
define a sequence of functions, $g_1,g_2,\dots, g_n$ from $(0,\infty)$ into
$[0,\omega_1),$ and a non-increasing 
sequence of ordinals $\gamma_1,\dots,\gamma_n$, then
 $f_{\epsilon,\mu}(t)$ will be $\sum_{i=1}^n\gamma_i+g_n(t)$ for some $n$
and all $t\leq \epsilon$.

First we assume
that $\mu (K^{d(\alpha)})\neq 0$ for only finitely
many $\alpha.$
Let $\alpha_1>\alpha_2>\dots>\alpha_k$ be the finite sequence of ordinals
such that $\lambda_i=\mu(K^{d(\alpha_i)})>0$ for each $i$ and
$\mu(K)=\sum_{i=1}^k \lambda_i$ and define $g_1(t)=\alpha_i$ if
$\sum_{j=1}^{i-1} \lambda_j<t \leq \sum_{j=1}^{i}\lambda_j,$ and $g_1(t)=0$
for $t>\sum_{i=1}^k \lambda_i.$ 

Before giving a formal description of the inductive procedure, let us
consider the following intuitive idea for a constructive approach to finding
the $\epsilon$-area. Notice that  the graph of $g_1$ is
decreasing. We would like to take the largest ordinal $\beta$ such that
$g_1(\epsilon)\geq \beta$, i.e., $g_1(\epsilon)$, let
$\gamma_1=\beta$ and define a new function $g_2$ as the decreasing
rearrangement of $g_1-1_{(0,\epsilon]}\gamma_1.$ The rectangle of width
$\epsilon$ and height $\gamma_1$ is our first approximation to the area
under $g_1$ and the region under $g_2$ is the remainder. Next we would
apply the procedure
to $g_2$ to get a new ordinal $\gamma_2=g_2(\epsilon)$ and let
$g_3$ be the decreasing rearrangement of $g_2-1_{(0,\epsilon]}\gamma_2.$ 
Proceeding inductively, we would find $(g_i)$ and $(\gamma_i)$. Notice that
$\gamma_i\geq \gamma_{i+1}$ and for only finitely many $i$ can we have
equality. Thus at some stage $\gamma_n=0$ and the procedure produces
nothing new.

Because of  some features of ordinal addition, it turns out that this
procedure may produce a little smaller function than we would like. To
avoid this it is necessary to require that $\gamma_i=\omega^{\beta_i}$ for
some $\beta_i$ and thus $\gamma_i$ may be strictly smaller than
$g_i(\epsilon).$ In the formal procedure below we will also describe in
detail a method for obtaining the decreasing rearrangement which will
allow us to extract some additional information for use later.
The main step in the
procedure is contained in the following lemma. Recall that if $\gamma$ and
$\beta$ are ordinals such that $\beta<\gamma$ then $\gamma-\beta$ is the
ordinal $\rho$ such that $\beta+\rho=\gamma.$ (See \cite{Haus}, page 74.)
In the statement of the lemma
and below $\lambda$ denotes Lebesgue measure.

\begin{lemma}
\label{rearrange-step}
Suppose that $g$ and $h$ are left continuous
non-increasing functions from $(0,\infty)$ 
into $[0,\omega_1)$ such that there exists $A<\infty$ with $g(t)=0=h(t)$
for all $t>A$,  
$g(t)\leq h(t)$ for all $t$, and the range
of each is a finite set of ordinals. Let $I=(a,b]$ be an interval on
which $g$ and $h$ are constant and let $\gamma \leq g(t)$ for $t\in I.$ Then
if $G$ and $H$ are the non-increasing left-continuous
rearrangements of $g-\gamma 1_I$ and
$h-\gamma 1_I$, respectively, then $G(t)\leq H(t)$ for all $t$ and
$\lambda(\{t:g(t)+1\leq h(t)\})\leq
\lambda(\{t:G(t)+1\leq H(t)\}).$
\end{lemma}
\begin{proof}
Let $s=\sup\{t:g(b)-\gamma<g(t)\}$ and
$r=\sup\{t:h(b)-\gamma<h(t)\}$.
Because $g$ and $h$ are non-increasing, 
\[G(t)=\begin{cases} g(t)& \text{ if  $t\leq a$ or $t>s$} \\
g(t+(b-a))& \text{ if $a<t\leq s-(b-a)$} \\
g(t-(s-b))-\gamma& \text{ if $s-(b-a)<t\leq s$}
\end{cases}
\]
and
\[H(t)=\begin{cases} h(t)& \text{ if  $t\leq a$ or $t>r$} \\
h(t+(b-a))& \text{ if $a<t\leq r-(b-a)$} \\
h(t-(r-b))-\gamma& \text{ if $r-(b-a)<t\leq r.$}
\end{cases}
\]
Observe that $G(t)\leq H(t)$ for all $t\leq p=\min (r,s)-(b-a)$ and
\[\lambda(\{t\leq p:g(t)+1\leq h(t)\})=\lambda(\{t\leq p:G(t)+1\leq H(t)\}).
\]
Similarly, if $q=\max(r,s)$, 
\[\lambda(\{t> q:g(t)+1\leq
h(t)\})=\lambda(\{t>q:G(t)+1\leq H(t)\}).\]
To see that $G(t) \leq H(t)$ for $q\geq t>p$ we note that if we do the same
rearrangement of $g-\gamma 1_I$ as for obtaining $H$ from $h-\gamma 1_I$,
we get
\[G_1(t)=\begin{cases} g(t)& \text{ if  $t\leq a$ or $t>r$} \\
g(t+(b-a))& \text{ if $a<t\leq r-(b-a)$} \\
g(t-(r-b))-\gamma& \text{ if $r-(b-a)<t\leq r.$}
\end{cases}
\]
and clearly all of the conclusions hold for $G_1$ in place of $G$. Now if
$G_1$ is not non-increasing, we
have two cases to consider. 
If $r<s$, then $G(t)=G_1(t-(s-r)) \leq G_1(t)) \leq H(t)$
for $s-(b-a)<t\leq s$ and
$G(t)=G_1(t+(b-a))\leq G_1(t)
\leq H(t) $ for $r-(b-a)<t\leq s-(b-a).$ Also 
\begin{multline*}
\{t:G(t)+1\leq H(t)\}\\
\supset\{t:G_1(t)+1\leq H(t),r-(b-a)<t\leq s-(b-a)\}\cup
(s-(b-a),s]. 
\end{multline*}
Thus the conclusion holds in this case.
 
If $r>s$, then $G(t)=G_1(t+(r-s))\leq H(t+(r-s))\leq H(t)$
for $s-(b-a)<t\leq s,$ and
$G(t)=G_1(t-(b-a))\leq G_1(t)\leq H(t)$ for $s<t\leq r.$ In this case 
\begin{multline*}
\{t:G(t)+1\leq H(t)\}\\
\supset\{t:G_1(t+(r-s))+1\leq H(t+(r-s)),s-(b-a)<t\leq s\}\cup(s,r]  
\end{multline*}
and the conclusion holds here too.
\end{proof}

The next lemma follows from a finite number of applications of
Lemma~\ref{rearrange-step}.

\begin{lemma}
\label{rearrange}
Suppose that $g$ and $h$ are left continuous
non-increasing functions from $(0,\infty)$ 
into $[0,\omega_1)$ such that $g(t)=0=h(t)$ for all $t>A$ for some $A,$ 
$g(t)\leq h(t)$ for all $t$, and the range
of each is a finite set of ordinals. Let $\epsilon>0$ and $\gamma>0$ such that
$\gamma \leq g(\epsilon)$.
Then
if $G$ and $H$ are the non-increasing rearrangements of $g-\gamma 1_{(0,\epsilon]}$ and
$h-\gamma 1_{(0,\epsilon]}$, respectively, then $G(t)\leq H(t)$ for all $t$ and
\[\lambda(\{t:g(t)+1\leq h(t)\})\leq
\lambda(\{t:G(t)+1\leq H(t)\}).
\]
\end{lemma}
\begin{proof} There are a finite number of disjoint, left-open, right closed
intervals $I_j, j=1,2,\dots,J$,
such that $1_{(0,\epsilon]}$, $g$ and $h$ are constant on each, and
$\cup_{j=1}^J I_j=\text{supp }h.$
We may assume that the intervals are ordered so that if $j_1<j_2$, $s\in
I_{j_1}$, and $t\in I_{j_2}$, then $s>t.$ There is some smallest index $j_0$
such that $I_{j_0}\subset (0,\epsilon]$. Applying Lemma~\ref{rearrange-step} to
$I_{j_0}$, $g$ and $h$, we get rearrangements $g^{(1)}$ and $h^{(1)}$ of
$g-\gamma 1_{I_{j_0}}$ and $h-\gamma 1_{I_{j_0}}$, respectively. Next we
repeat the process with $I_{j_0+1}$, $g^{(1)}$ and $h^{(1)}$ to obtain
$g^{(2)}$ and $h^{(2)}$. Clearly, this process produces the required
non-increasing rearrangements of $g-\gamma 1_{(0,\epsilon]}$ and $h-\gamma
1_{(0,\epsilon]}$ at stage $J-j_0+1$. Because $g^{(j)}\leq h^{(j)}$ for
and 
\[ \lambda(\{t:g^{(j)}(t)+1\leq h^{(j)}(t)\})\leq
\lambda(\{t:g^{(j+1)}(t)+1\leq h^{(j+1)}(t)\}).
\]
for each $j$, the required properties follow immediately
\end{proof}

\begin{remark} \label{rem1}
Notice that if $h$ is non-increasing as in Lemma \ref{rearrange} and $\rho$
is an ordinal such that $\rho\cdot\omega<h(\epsilon)$, then $h-\rho
1_{(0,\epsilon]}=h$. Thus if too small an ordinal is chosen, there is no
effect.
\end{remark}

The next proposition will enable us to define the $\epsilon$-distribution.
Below we use summations of ordinals with the understanding that
$\sum_{i=1}^n \gamma_i= \gamma_1+\gamma_2+\dots+\gamma_n$ in that order.

\begin{proposition}
\label{dist}
Let $\epsilon>0$ and let $g_0:(0,\infty)\rightarrow [0,\omega_1)$ be a
left continuous,
non-increasing function with range a finite set
such that for some $t_0<\infty,$ $g_0(t)=0$ for all $t>t_0.$
Then there exists a finite sequence of left continuous, non-increasing
functions $(g_i)_{i=1}^n$ from $(0,\infty)$ into
$[0,\omega_1)$ and a non-increasing sequence of ordinals
$(\gamma_i)_{i=0}^{n-1}$ such that for each $i<n$ and $\alpha<\omega_1$,
\[\lambda(\{t:g_{i+1}(t)=\alpha\})=\lambda(\{t:g_i(t)-\gamma_i
\text{\large 1}_{(0,\epsilon]}(t) = \alpha\}),\]
i.e., $g_{i+1}$ is a decreasing rearrangement of $g_i-\gamma_i
1_{(0,\epsilon]},$ $\gamma_i=\omega^{\beta_i}$ for some $\beta_i,$
and $g_n(t)=0,$ for all $t\geq\epsilon$.

Moreover, if $g_0$ and $h_0$ are two non-increasing
functions as above, $g_0(t)\leq h_0(t)$ for all $t$,
and $(g_i)_{i=1}^n$, $(\gamma_i)_{i=1}^{n-1}$, and
$(h_i)_{i=1}^{m}$, $(\eta_i)_{i=1}^{m-1}$, are the corresponding sequences
of functions and ordinals produced, then $\sum_{i=1}^{n-1} \gamma_i
+g_n(t)\leq \sum_{i=1}^{m-1} \eta_i+h_m(t),$ for all $t\leq \epsilon.$
Further, if $\lambda\{t:g_0(t)+1\leq  h_0(t)\} \geq \epsilon$,
then $\sum_{i=1}^{n-1} \gamma_i+g_n(\epsilon)+1\leq h_m(\epsilon).$
\end{proposition}

\begin{proof}
The proof proceeds by constructing inductively the sequence $(g_i)$. In
order to prove the moreover assertion we will work with $h_0$ at the same
time
and produce the corresponding sequence $(h_i)$.

Suppose that we have $g_i$ and $h_i$,
$1,2,\dots k,$ such that $g_i\leq h_i$ for each $i$. If
$g_k(\epsilon)=0,$ the
construction of the sequence $(g_i)$ is complete.
If not let $\beta_k$ be the largest ordinal $\beta$
such that $\omega^\beta\leq g_k(\epsilon)$.
Let $g=g_k$, $h=h_k$, $\gamma=\omega^{\beta_k}$, and $I=(0,\epsilon].$
Applying Lemma~\ref{rearrange}
we let $g_{k+1}=G$ and $h_{k+1}=H$ be the decreasing
rearrangements of $g_k-\gamma_k 1_I$ and $h_k-\gamma_k 1_I$ such that $G\leq
H$. Moreover, $\lambda(\{t:g_k(t)+1\leq h_k(t)\})\leq
\lambda(\{t:g_{k+1}(t)+1\leq h_{k+1}(t)\})$.

Notice that if $h_k(\epsilon)\geq \gamma_k\cdot \omega,$ $h_k-\gamma_k
1_I=h_k.$
Thus if this occurs for some $k$, $h_k=h_i$ for all $i,$ $k\leq i\leq n$, and
\[\sum_j^{k-1} \eta_j+ h_k(\epsilon)> \sum_j^{i-1} \gamma_j+
g_i(\epsilon)+1
\]
for each $i$.
If  $h_k(\epsilon)< \gamma_k\cdot \omega,$ for all $k\leq n-1,$ then each step
of the construction of $(g_i)$ is also a step in the construction of $(h_i)$
with $\eta_k=\gamma_k$. Clearly,
\[ \sum_{j=n}^{i-1}\eta_j+h_i \geq g_n\]
for $i=n,n+1,\dots,m$.
This completes the proof of all of the conclusions except for
the final assertion in the case $h_k(\epsilon)<\gamma_k\cdot \omega$.

Because $\lambda(\{t:g_n(t)+1\leq h_n(t)\})\geq \epsilon$, at step $n$
either $h_n(t)=0$ for all $t>\epsilon$ and $h_n(t)\geq
g_n(t)+1$ for all $t$, $0\leq t \leq \epsilon,$ or $h_n(t)>0$ for some
$t>\epsilon.$ The first case satisfies the conclusion of the proposition.
In the second case observe that for each $i$, 
$\sum_{j=1}^i \gamma_j +h_{i+1}(\epsilon)\geq
\sum_{j=1}^{i-1}\gamma_j+h_i(\epsilon).$
Because $h_n(t)>0$ for some $t>\epsilon$, it follows that there is a
largest $\eta_n=
\omega^{\beta_n}>0$ such that $\eta_n\leq h_n(\epsilon)$. 
Because $\eta_n>g_n(\epsilon)$, the proof is complete.
\end{proof}

We now introduce terminology for
some of the ingredients of Proposition~\ref{dist} and its proof.
\begin{definition}\label{e-dist}
Suppose $g$ is a non-increasing left-continuous function from $(0,\infty)$
into $(0,\omega_1)$ and $\epsilon>0.$ If $\gamma\leq
g(\epsilon)$ and $f$ is the decreasing rearrangement of
$g-\gamma 1_{(0,\epsilon]}$ then $h=\gamma 1_{(0,\epsilon]}+f$ will be said
to be an {\it $\epsilon$-compression} of $g$ (by $\gamma$).

Let
\begin{align*}
C(g,\epsilon)= 
\sup &\{H(\epsilon):\text{ there exist non-increasing left-}\\
&\text{continuous} \text{ simple}
\text{ functions } (h_i)_{i=1}^n, h_1\leq g, \\
&h_{i+1}\text{ is an $\epsilon$-compression of }h_{i},
H=h_n\}
\end{align*}
For a positive finite measure $\mu$ on $K$ let $g(t)=\sup\{\alpha:
\mu(K^{(\alpha)})\geq t\}$ for all $t>0$
and define $C(\mu,\epsilon)=C(g,\epsilon).$
(We let the supremum of an empty set of ordinals be $0$.)
We will call $g$ the {\it derived height} of $\mu$ and $C(g,\epsilon)$ the
\textit{$\epsilon$-area} under $g$.
\end{definition}

It is not hard to see that the procedure used in the proof of Proposition
\ref{dist} will produce the value of $C(g,\epsilon)$ if $g$ is simple. In
that case with $g_j$ and $\gamma_j$ as in the proof we let
$h_i=\sum_{j=1}^{i-1} \gamma_j 1_{(0,\epsilon]}+g_i$. It is important
in achieving the supremum that for each
$j$, $\gamma_j$ is of the form $\omega^{\beta_j}$. This avoids lowering
the sum by taking the wrong order, e.g., $\omega^2+1$ and $\omega$ sum
(in that order) to $\omega^2+\omega$ but $\omega^2$, $\omega$, and $1$ sum to
$\omega^2+\omega+1$.

Observe that for a measure $\mu$ as
in Definition~\ref{e-dist}, if for some $t$, $g(t)=\alpha$,
$\mu(K^{(\alpha)})\geq t$. Also if $(t_n)$ is an increasing sequence of
positive numbers with limit $t$
and $g(t_n)=\alpha_n$ for each $n$, $(\alpha_n)$ must eventually be constant.
Thus $\mu(\cup
K^{(\alpha_n)})=\lim \mu(K^{(\alpha_n)})\geq\lim t_n$
and
$g(t)=\lim g(t_n),$ i.e., $g$ is left-continuous.

Also notice that if $g$ and $h$ are non-increasing functions as in the
statement of the proposition but not necessarily simple, then the final
conclusion of
Proposition~\ref{dist} still holds, i.e, $C(g,\epsilon)+1\leq C(h,\epsilon)$.
Indeed, if $g_1$ is a simple function with $g_1\leq g$ and $A$ is the set
where $h(t)\geq g(t)+1$, then $g_1+1_A \leq h$. It follows easily that there
is a non-increasing simple function $h_1$ such that $g_1+1_A \leq h_1\leq
h$. Thus $C(g_1,\epsilon)+1\leq C(h_1,\epsilon)\leq C(h,\epsilon)$. Taking
the supremum over all such $g_1$ gives the result.

\begin{example}
If we return to our previous example
\[g(t)=\omega \text{\large 1}_{(0,1/4]}+\text{\large 1}_{(1/4,1]} \]
and let $\epsilon =1/2,$ then $C(g,1/2)=2$ because $h_1=g$ and
\[ H=h_2=\omega \text{\large 1}_{(0,1/4]}+2\text{\large 1}_{(1/4,1/2]}.\] If $\epsilon=1/4$ then
$C(g,1/4)=\omega +3$. Indeed, $h_1=g$,
\begin{align*}
h_2&=(\omega+1)\text{\large 1}_{(0,1/4]}+\text{\large 1}_{(1/4,3/4]}\\
h_3&=(\omega+2)\text{\large 1}_{(0,1/4]}+\text{\large 1}_{(1/4,1/2]}\\
H=h_4&=(\omega+3)\text{\large 1}_{(0,1/4]}.
\end{align*}
\end{example}

\begin{remark}
The definition of $\epsilon$-area can be adapted to accommodate different
values of $\epsilon$  as in the definition of summable Szlenk index, \cite{GKL}
or \cite{KOS}, but one must use the differences instead of the
$\epsilon$-compressions. Thus one would begin with $g$ and $\omega^\gamma=1$
and let $g_1$ be the decreasing rearrangement of $g-1_{(0,\epsilon_1]}$,
$g_2$ be the decreasing rearrangement of $g_1-1_{(0,\epsilon_2]}$, etc.
The $(\epsilon_1,\epsilon_2,\dots,\epsilon_n)$-area is zero if
$g_i-1_{(0,\epsilon_{i+1}]}$ is not non-negative for some $i$. This notion
of summable Szlenk index seems to be the same as saying that there is a
constant $K$ such that for every $\epsilon>0$, the $\epsilon$ area or
equivalently the Szlenk index is at most $[K/\epsilon]+1$, where $[\cdot]$
denotes the greatest integer. (See \cite{KOS} where this latter property
is called proportional index.)
\end{remark}

Our next task is to show that there is a relation between the derivation on
$K$ and a ``Szlenk'' derivation on the probability measures on $K$. Below
the $\text{weak}^*$-topology on the probability measures on $K$ is that
inherited from $C(K)^*.$

\begin{definition}\label{M-derive}
Suppose that $M$ is a set of probability measures on $K$ and let
$\delta,\epsilon>0.$ Define $M(\epsilon,\delta)^{(0)}=M.$
For each $\alpha<\omega_1$, define 
\begin{align*}
M(\epsilon,\delta)^{(\alpha)}=&\{\mu:\text{ there exists
}(\mu_n)_{n=1}^\infty \subset M(\epsilon,\delta)^{(\alpha)}
\text{ and a sequence} \\
&\text{of closed subsets $(A_n)_{n=1}^\infty$
 of $K$ such that, }
\mu_n(A_n)\geq \epsilon \\
&\text{for all $n$, }\text{w}^*\lim \mu_n=\mu,\, d(A_n,A_m)\geq
\delta \text{ for all }n\neq m\}.
\end{align*}
If $\beta$ is a limit ordinal, define $M(\epsilon,\delta)^{(\beta)}
=\cap_{\alpha<\beta}
M(\epsilon,\delta)^{(\alpha)}.$
\end{definition}

Notice that the definition is at least superficially
more restrictive than that of the $\epsilon$-Szlenk subsets
of $M$ in that from the Szlenk index definition we
would only have disjointness of the sets $(A_n)$ not separation by $\delta$.
Indeed, if $(\mu_n)$ is a $\text{w}^*$ convergent sequence of probability
measures and $(f_n)$ is a weakly null sequence of
(without loss of generality) positive continuous
functions such that
$\int f_n\,d\mu_n\geq \epsilon$,  then given $\epsilon'<\epsilon$,
for a sufficiently small $\rho>0$,
we can let $A'_n=\{k:f_n(k)>\rho\}$ for each $n$ and by passing to a
subsequence if necessary, let $A_n=A'_n\setminus \cup_{k<n}A_k$,
to obtain disjoint sets
such that $\int_{A_n} f_n \,d\mu_n>\epsilon',$ for all $n$. Essentially
this is the same as saying that the Szlenk definition detects the
non-uniform absolute continuity of a set of measures. (See \cite{A1} and
the proof of Corollary~\ref{cor-ord}.)

\begin{proposition} \label{Cep}
Let $\epsilon, \delta>0.$ If $M$ is a subset of the probability measures on 
a compact set $K$ with $\delta$-system of derived sets
$\{K^{(\alpha)}:\alpha<\omega_1\}$
and $\mu\in M(\epsilon,\delta)^{(\alpha)}$,
then for every $\epsilon'<\epsilon$,
$C(\mu,\epsilon')\geq \alpha.$
\end{proposition}

\begin{proof}
The proof is by induction on $\alpha$. The main step is to prove the
following.

\medskip
\noindent{\bf Claim:}  If $(\mu_n)_{n=1}^\infty \subset M$ with
$C(\mu_n,\epsilon')\geq\alpha$ for each $n$, $\text{w}^*\lim \mu_n=\mu$, 
 and $(A_n)_{n=1}^\infty$ is a sequence of closed subsets
 of $K$ such that $\mu_n(A_n)\geq \epsilon$ and
 $d(A_n,A_m)\geq \delta$  for all $n\neq m$,
then $C(\mu,\epsilon')\geq \alpha+1.$

\medskip
Because the sets $K^{(\beta)}$ are closed for each $\beta$ and $\mu_n\geq
0$, 
\[\limsup
\mu_n(K^{(\beta)}) \leq \mu(K^{(\beta)}) 
\]
for all $\beta.$ Therefore if $h_n$
is the derived height of $\mu_n$ for each $n$ and $h$ is the derived height of
$\mu$, then $h(t)\geq \limsup h_n(t)$ for all $t$. Given $\rho>0$ we can find
$\alpha_1<\alpha_2<\dots<\alpha_k$ such that $\sum_{i=1}^k
\mu(K^{d(\alpha_i)})>1-\rho.$ (Let $\alpha_0=0$.) Now consider 
$\delta_{i+1}=\limsup \mu_n((K^{(\alpha_i)}\setminus
K^{(\alpha_{i+1})})\cap A_n)$. By passing to a subsequence we may assume
that this limit exists for each $i$ and so does $\lim \mu_n(K^{(\alpha_i)}).$
If $k_n\in (K^{(\alpha_i)}\setminus
K^{(\alpha_{i+1})})\cap A_n$, we know that any limit point of $(k_n)$ is in
$K^{(\alpha_i+1)}.$
Therefore, if $g$ is a continuous function such that
$1_{K^{(\alpha_i+1)}}\leq g \leq
1$ and $\rho'>0$,
then $(\mu_n,g)\geq \mu_n(K^{(\alpha_{i+1})})+\mu_n(K^{(\alpha_i)}\setminus
K^{(\alpha_{i+1})})\cap A_n)-\rho'$ for $n$ sufficiently large.
Consequently, 
\begin{multline*}
\mu(K^{(\alpha_{i}+1)}\setminus
K^{(\alpha_{i+1})})+\mu(K^{(\alpha_{i+1})})=\mu(K^{(\alpha_i)})\\
\geq \limsup
\mu_n(K^{(\alpha_{i+1})})+\delta_{i+1}.
\end{multline*}
Rearranging, we get
\[\mu(K^{(\alpha_{i+1})})-\limsup \mu_n(K^{(\alpha_{i+1})})\geq
\delta_{i+1}-\mu(K^{(\alpha_{i}+1)}\setminus
K^{(\alpha_{i+1})}).
\]

Because 
\begin{align*}
\sum_{i=1}^{k} \delta_i &=\sum_{i=1}^k \limsup
\mu_n((K^{(\alpha_{i-1})}\setminus K^{(\alpha_{i})})\cap A_n)\\
&\geq \limsup \mu_n(A_n \cap (K \setminus K^{(\alpha_k)}))\\
&\geq \limsup \mu_n(A_n\cap K)-\mu_n(A_n\cap K^{(\alpha_k)})\\
&\geq \epsilon-\mu(K^{(\alpha_k+1)})\geq \epsilon-\rho,
\end{align*}
\begin{multline*}
\sum_{i=0}^{k-1} \mu(K^{(\alpha_{i+1})})-\lim \mu_n(K^{(\alpha_{i+1})})\\
\geq
\sum_{i=0}^{k-1}
\delta_{i+1}-\mu(K^{(\alpha_{i}+1)}\setminus
K^{(\alpha_{i+1})})\geq \epsilon-2\rho.
\end{multline*}
Clearly $h(t)\geq \limsup h_n(t)+1$
for $t$ such that $\lim \mu_n(K^{(\alpha_i)})<t\leq
\mu(K^{(\alpha_i)})$.
Thus 
\[\lambda(\{t:h(t)\geq \limsup h_n(t)+1
\})\geq \epsilon-2\rho,\]
 for every $\rho>0.$ It follows there is some $n$
such that $\lambda(\{t:h_n(t)+1\leq h(t)\})>\epsilon'.$ Proposition~\ref{dist}
implies
that $C(h_n,\epsilon')+1\leq C(h,\epsilon'),$ proving the Claim.

The Claim proves the induction step. Indeed, if $\mu \in
M(\epsilon,\delta)^{(\alpha+1)}$ then
there is a sequence $(\mu_n)\subset M(\epsilon,\delta)^{(\alpha)}$
with $\text{w}^*$-limit $\mu$ as
in the Claim. By the inductive assumption $C(\mu_n,\epsilon')\geq \alpha$
and thus the Claim gives $C(\mu,\epsilon')\geq \alpha+1.$

If $\alpha$ is a limit ordinal, let $(\alpha_n)$ be a sequence of ordinals
converging to $\alpha.$ If $\mu\in M(\epsilon,\delta)^{(\alpha)}$,
then $\mu \in
M(\epsilon,\delta)^{(\alpha_n)}$ for all $n$. By the induction
hypothesis $C(\mu,\epsilon')\geq
\alpha_n$ for all $n$. Therefore $C(\mu,\epsilon')\geq \alpha.$
\end{proof}

The next result is known, e.g., \cite{S}, but the apparatus we have
constructed gives an easy proof.

\begin{corollary} \label{cor-ord}
The $\epsilon$-Szlenk index of the unit ball of
$C(\omega^{\omega^\gamma\cdot k})^*$ is 
$\omega^\gamma [k/\epsilon]+1$.
\end{corollary}
\begin{proof}
We take the $\delta$-system of derived sets
to be the usual topological derived sets of $[1,\omega^{\omega^\gamma\cdot k}] \cup
-[1,\omega^{\omega^\gamma\cdot k}]$, (the disjoint union of two copies of
$[1,\omega^{\omega^\gamma\cdot k}]$), the
metric to be the discrete metric $d(x,y)=1,$ for $x\neq y$, and $\delta=1$.
If $\mu$ is any probability measure on  $[1,\omega^{\omega^\gamma\cdot k}]\cup
-[1,\omega^{\omega^\gamma\cdot k}]$, the
derived height can be at most $\omega^\gamma\cdot k$ at each point.
Thus $C(\mu,\epsilon)\leq \omega^\gamma [k/\epsilon]$.

Now consider the definition of the Szlenk subsets of the ball of
$C([1,\omega^{\omega^\gamma\cdot k}])^*$, $P_\alpha(\epsilon')$. If
$\mu\in P_{\alpha+1}(\epsilon')$ then there is a sequence of measures $(\mu_n)$
which converge $\text{w}^*$ to $\mu$ and a sequence of norm one continuous functions
$(f_n)$ converging pointwise to $0$ such that $\lim (\mu_n,f_n)\geq
\epsilon'.$ Let $\epsilon''<\epsilon'$.
It follows that there are disjoint sets $(A_n)_{n\in K}$ for
some infinite set $K\subset \mathbb N$ such that $|\mu_n|(A_n)\geq \epsilon
''$. Except for the absolute values this is precisely
the condition in Definition~\ref{e-dist}. We can eliminate the absolute
values by considering measures on $[1,\omega^{omega^\gamma\cdot k}]
\cup-[1,\omega^{\omega^\gamma\cdot k}]$. (Replace $\mu$ by $\mu'$ where
$\mu'(A)=\mu^+(A\cap [1,\omega^{\omega^\gamma \cdot k}])
+\mu^-(-(A\cap-[1,\omega^{\omega^\gamma\cdot k}])).$) Thus if $\mu
\in P_\alpha(\epsilon'),$ then $\mu \in M(1,\epsilon'')^{(\alpha)},$ 
where $M$ is the set of probability measures on $\pm[1,\omega^{\omega^\gamma\cdot
k}]$.
Thus to compute the Szlenk index 
we may apply Proposition \ref{Cep} to get that $\mu \in
P_{\alpha}(\epsilon')$ implies that $C(\mu,\epsilon'')\geq \alpha.$
Therefore $\alpha\leq \omega^\gamma [k/\epsilon'']$ for every
$\epsilon''<\epsilon',$  and the $\epsilon'$ Szlenk index is at most
$\omega^\gamma [k/\epsilon']+1.$
It is easy to see that $\delta_{\omega^{\omega^\gamma\cdot k}}\in P_{\omega^\gamma
[k/\epsilon']},$ completing the proof.
\end{proof}

\section{The Szlenk Index of the Bourgain-Delbaen Space}
\label{szcomp}
The proof of Proposition \ref{quo}
suggests the following approach to representing
(non-uniquely) the $\text{w}^*$-closure of the basis $\{e_k^*:k\in \mathbb N\}$.
Recall that a tree is a partially ordered set $(T,\leq)$ such that each
initial segment, $\{y:y\leq x\}$ for $x\in T,$ is well-ordered and finite.
Let $T=\cup_{n=0}^\infty  \{0,1\}^n$, the rooted binary tree (ordered by
extension) with root the
empty tuple, $()$,
 and let 
\[W= (\{0,a,-a,b,-b,1\}\times\{\omega\cdot
m+k:m,k\in \mathbb N\cup \{0\}\})\cup\{\infty\},\]
the one-point compactification of
$\{0,a,-a,b,-b\}\times [1,\omega^2)$. (In this topology any sequence in $W$
of the form $(c_i,\omega\cdot m_i+k_i)$ with $\lim m_i=\infty,$ has limit
$\infty.$) Let $K$ be the space of all functions
from $T$ into $W$ in the topology of pointwise convergence. We have that $K$ is
compact by the Tychonoff theorem. Each basis vector $e_k^*$ in $X^*$ can be
associated to a point $g_k$ in $K$ in the following way. 

Let $g_k(())=(1,k)$ and if 
$g_k(\delta_1,\delta_2,\dots,\delta_n)$ has been defined to be 
$(c,\omega\cdot m+\ell)$ and
$\phi(\ell)=(\sigma',i,m',\sigma'',j)$, let
\[g_k(\delta_1,\delta_2,\dots,\delta_{n+1})=\begin{cases}
(\sigma'a,\omega\cdot
m+i) &\text{ if $\delta_{n+1}=0$}, \\
(\sigma'' b, \omega \cdot \max(m,m')
+ j) &\text{ if $\delta_{n+1}=1$}.
\end{cases}
\]
If $\ell\leq 2$, \[g_k(\delta_1,\delta_2,\dots,\delta_{n+1})=(0,\omega\cdot
m).\]
 Define $\theta(g_k)=e_k^*,$ for all $k$.

We need some notation to conveniently refer to the pieces of $W$.
For $(c,\omega\cdot m+j)$ we define three functions which extract
the essential parts: $V(c,\omega\cdot m+j)=c$, $Q(c,\omega\cdot
m+j)=m$, and $R(c,\omega\cdot m+j)=j$. For a node
$\mathcal N$ of the binary tree of length $L(\mathcal N)=n$
and $t<n$ define the $t$th
truncation by $I(\mathcal N,t)=(\delta_1,\delta_2,\dots,\delta_t)$ if
$\mathcal N=(\delta_1,\delta_2,\dots,\delta_n)$.
We define the evaluation of an element $x$ of $X$ by an element $f\in K$
at a node $\mathcal N=(\delta_1,\delta_2,\dots,\delta_n)$ by
\[<f,\mathcal N,x>=\left (\prod_{j=1}^n
V\big (f(I(\mathcal N,j))\big )\right )
\left ((I-P_{Q(f(\mathcal N))}^*)
e_{R(f(\mathcal N))}^*\right ) x.\]
Now suppose that $x\in P_sX$ for some $s$ and $f$ is the preimage of $e_k^*$
for some $k$, i.e., $f=g_k$. For each node $\mathcal B=(\delta_i)$ of
$T$ there is smallest index $n=n(\mathcal B,f)$
such that $R(f(I(\mathcal B,n)))\leq d_s.$ (Of course every node with this
initial segment yields the same index.)

\begin{proposition} Let $f=g_k$, i.e., $\theta(f)=e_k^*$,
 for some $k$ and $x\in P_sX$
for some $s$. If $\{\mathcal N_i\}$ is a maximal
collection of incomparable nodes such that
$L(\mathcal N_i)\leq n(\mathcal B,f)$ for any branch $\mathcal B$ with
$\mathcal N_i$ as an initial segment.
Then the collection is finite and $e_k^*(x)=\sum_i (f,\mathcal N_i,x).$ 

\end{proposition}
\begin{proof}
Observe that if $\mathcal N$ is any node with $L(\mathcal
N)<n(\mathcal B,f)$, $f(\mathcal N)=(c,\omega\cdot m+j)$ and $\phi(j)=
(\sigma',r,m',\sigma'',q)$, then  $j>d_s$
and
\begin{multline}
\label{nodeval}(I-P_m^*)e_j^*(x) =\\
\sigma'a(I-P_m^*)e_r^*(x)+\sigma''b(I-P_{m'}^*)(I-P_m^*)e_q^*(x).
\end{multline}
Note that $(I-P_{m'}^*)(I-P_m^*)=I-P_{\max(m,m')}^*.$
If $(f,\mathcal N,x)=c(I-P_m^*)e_j^*(x)$, then
\begin{equation}\label{nodesplit}
\begin{split} (f,\mathcal
N,x)&=c(\sigma'a(I-P_m^*)e_r^*(x)+\sigma''b(I-P_{\max(m,m')}^*)e_q^*(x))\\
&=
(c\sigma'a)(I-P_m^*)e_r^*(x)+(c\sigma''b)(I-P_{\max(m,m')}^*)e_q^*(x)\\
&=(f,\mathcal N\mathbf{+}(0),x)+(f,\mathcal N\mathbf{+}(1),x),
\end{split}
\end{equation}
where $(\cdot)\mathbf{+}(\cdot\cdot)$ denotes the concatenation of
the tuples $(\cdot)$ and $(\cdot\cdot)$.
  Therefore we can
prove the formula by induction on the set of nodes as follows.
We enumerate  the nodes
of the binary tree so that all nodes of a given length are labeled before
any node of a longer length. Observe that the formula is obvious if we have
only the node $()$ since 
\[(f,(),x)=(I-P_0^*)e_k^*(x)=e_k^*(x).\]
If this is the
maximal collection, we are finished. If not, $()$ is the first node in the
enumeration and we replace it by the two node
collection $\{(0),(1)\}$. Formula \eqref{nodesplit}
immediately gives the result if this
is the collection of nodes. Otherwise we consider the next node in the
enumeration. If it is in the collection $\{\mathcal N_i\}$,
 we move on in the enumeration;
if not we apply the
formula \eqref{nodesplit} to
replace the node by the two nodes immediately below.
Note that because we began with $e_k^*$, with $k\leq d_r$ for some $r$,
 the integer coordinates of
$\phi(k)$ are smaller than $d_{r-1}$. Iterating, we see that
 there can be only finitely many
nodes in the collection $\{\mathcal N_i\}$.
Continuing in this way we eventually reach each node in the original
collection and the formula follows.
\end{proof}

Our next task is to show that if $\{g_k\}$ is the set of
representatives in $K$ of the basis elements $\{e_k^*\}$ defined above, then
the mapping $\theta$ described above extends to a continuous map from
$\overline{\{g_k\}}$ into $X^*.$ 

Before we proceed, let us note that because of the role of $m=Q(g_k(\mathcal
N))$, once there
is a node $\mathcal N_0$ in a branch that contains $0$, 
\[R(g_k(\mathcal M))\leq Q(g_k(\mathcal N_0))
\]
for all nodes $\mathcal M$ which are
descendants of $\mathcal N_0$.
Hence there can be only finitely many
nodes on the branch containing $\mathcal N_0$ at which $V$ is non-zero.

\begin{proposition}
The map $\theta$ extends to a continuous function from $\overline{\{g_k\}}$
into $\overline{\{e_k\}}$.
\end{proposition}
\begin{proof}
Suppose that $(g_k)_{k\in M}$ has limit $f$ in $K$. We have that 
$(g_k(\mathcal N))$ converges for each node $\mathcal N.$
$g_k(\mathcal N)=(c_k,\omega\cdot m_k+j_k)$ for each $k$. If $(m_k)$ is not
bounded then $\lim m_k=\omega$ and
limit of $(g_k(\mathcal N))$ is $\infty$. Assume that this is not the case.
Because for each $k$, $c_k$ and $m_k$ must
be one of a finite set of values it follows that $(c_k)$ and $(m_k)$ are
eventually constants $c$ and $m$, respectively.
If $(g_k(\mathcal N))$ is not eventually constant then
$\lim_{k\in M} j_k = \omega$ and the limit is $(c,\omega\cdot (m+1))$.
Therefore for each node we have three possible situations.
\begin{enumerate}
\item{} $(g_k(\mathcal N))$ converges to $\infty$.
\item{} $(g_k(\mathcal N))$ is eventually constant.
\item{} $(g_k(\mathcal N))$ converges to $(c,\omega\cdot (m+1)).$
\end{enumerate}
Consider in each case what happens on the nodes below.

In the first and second cases by \eqref{nodeval}
the same must be true for each node below $\mathcal N.$ In the third case
we must exam $(\phi(j_k))$ as in the proof of the Proposition~\ref{quo}. 
Observe that $(g_k,\mathcal N,x)=c_k(I-P_{m_k}^*)
e_{j_k}^*(x))$ for some constant $c_k$
and consider the same three cases. In the first case $(m_k)$
diverges to $\infty$ and therefore $\lim c_k(I-P_{m_k})x=0$ for every $x\in 
\cup_s P_s E_s.$ Consequently, $\text{w}^*\lim c_k(I-P_{m_k}^*)e_{j_k}^*=0.$

In the second case $((g_k,\mathcal N,x))$ is eventually constant and so is
$(c_k(I-P_{m_k}^*)e_{j_k}^*)$.

In the third case $(c_k)$ and $(m_k)$ are eventually constant and
consequently, so is $(f,\mathcal N\mathbf{+}(0),x)$.

To determine the limit of $\theta(g_k)$ we let $\{\mathcal N_i\}$ be the
sequence of nodes such that $f(\mathcal N_i)\neq \infty$, 
$R(f(\mathcal N_i))\neq 0$ and $R(f(I(\mathcal
N_i,L(\mathcal N_i)-1))=0.$ By definition this is a set of
incomparable nodes. Define $y^*(x)=\sum_i (f,\mathcal N_i,x).$ We claim
that $\text{w}^*\lim \theta(g_k)=y^*.$ Indeed the nodes we have described above
are precisely the nodes corresponding to the terms
that appear in the series representation for a
limit point of $\theta(g_k)$ determined in the proof of Proposition~\ref{quo}.

\end{proof}

Let $C=\overline{\{e_k^*\}}.$
\begin{proposition}\label{Sz-C}
For each $\epsilon>0$ the $\epsilon$-Szlenk index $\eta(\epsilon,C)$ is finite.
\end{proposition}
\begin{proof}
Fix $\epsilon>0.$
Find $N$ such that 
\[\sum_{i=N}^\infty a^i \sup_{m,k}\|(I-P_m)e_k^*\|
<\epsilon/4.\]
Suppose that $(x_k^*)$ is an $\epsilon$-separated sequence in $C$ and
\[x_k^*=\sum_{j=1}^\infty c_{k,j}(I-P_{m_{k,j}})e_{i_{k,j}}^*\]
 and $|c_{k,j}|
\leq a^j$ for all $j$. Then $y_k=\sum_{j=1}^{N-1}
c_{k,j}(I-P_{m_{k,j}})e_{i_{k,j}}^*$, $k=1,2,\dots,$
is an $\epsilon/2$-separated sequence.
Let $(z_k)$ be the sequence of preimages of $(x_k)$ corresponding to the
series representation above. ($\theta(z_k)=x_k$ for
all $k$.) Because the $(y_k)$ is $\epsilon/2$ separated it follows that
the sequence of restrictions $(z_k|_{\{\mathcal N:L(\mathcal N)<N\}})$ is
distinct. Now observe that the set of maps from  a finite set $G$ into $W$
in the topology of pointwise convergence is a metric space homeomorphic to
$[1,\omega^{2\cdot \text{ card }G}(6\cdot\text{ card }G)]$. Therefore the
$\epsilon$-Szlenk index is at most $2^{N+1}+1.$
\end{proof}

\begin{corollary} For each $\epsilon>0$, $\eta(\epsilon,B_{X^*})<\omega.$
Consequently, $C(\omega^\omega)$ is not isomorphic to a quotient of $X$.
\end{corollary}

\begin{proof}
It is sufficient to consider $D$=$\overline{\text{co }\pm\{e_k^*:k\in \mathbb
N\}}^{\text{w}^*}$ in place of $B_{X^*}.$ By the Choquet theorem we can
associate each element $x^*$ of $D$ to some probability  
measure  $\mu_{x^*}$ on
\[
C=\overline{\pm\{e_k^*:k\in \mathbb N\}}^{\text{w}^*}.
\]
Observe that if
$(x_n^*)$ is a $\text{w}^*$-convergent sequence in $D$ and $(x_n)$ is a
weakly null sequence in the unit ball of $X$ such that $\lim
x_n^*(x_n)\geq \epsilon_1$, then there exist an infinite subset $L$ of
$\mathbb N$ and $\epsilon/4$ norm separated subsets $(A_n)_{n\in L}$ of $C$
such that 
$\mu_{x_n^*}(A_n)\geq \epsilon_1/2$ for all $n\in L$.

Now we consider the modified Szlenk subsets of $C$,
$\{P_\alpha(\epsilon_1/4,C):\alpha<\omega_1\}$,
as the $\delta$-system
of derived sets with $\delta=\epsilon_1/4$ and  let 
\[M=\{\mu: \mu\text{ is a probability measure
representing some }x^*\in D\} ,
\]
$\epsilon=\epsilon_1/2$
 and $\delta=\epsilon_1/4$. By Proposition
\ref{Cep} if $\mu\in M(\epsilon,\delta)^{(\alpha)}$
then $C(\mu, \epsilon_1/4)\geq \alpha$.
However $C(\mu,\epsilon_1/4)$ must be finite by Proposition \ref{Sz-C}.
\end{proof}

\begin{remark} From the proofs of Proposition \ref{Sz-C} and the corollary,
the $\epsilon$-Szlenk index can be estimated from above. Recently Haydon
\cite{Haydon} has shown that the Bourgain Delbaen spaces are hereditarily
$\ell_p$ for some $p$ which depends on $a$ and $b$. From this one can get a
lower estimate on the Szlenk index. Using results in \cite{GKL,GKL1} it follows
that these spaces are not uniformly homeomorphic to $c_0$.
\end{remark}

\begin{remark} After reading an earlier version of this paper
I. Gasparis communicated to us another method of showing that the
$\epsilon$-Szlenk index of the Bourgain-Delbaen space is finite without
determining the behavior of the index. With his
permission we include a sketch of the argument here.

$\eta(\epsilon,B_{X^*})\geq \omega$ for some $\epsilon>0$ is equivalent to the
statement that $C(\omega^\omega)$ is a quotient of $X$. (See \cite{AB}.) It
is well-known that
$C(\omega^\omega)$ has $\ell_1$ as a spreading model of a weakly null
sequence. If $C(\omega^\omega)$ is a quotient of $X$, then $X$ also
has a weakly
null sequence with spreading model $\ell_1.$ This would imply that the
basis of $X$ has blocks that are equivalent to the basis of
$\ell_1^n$ for all $n$. However the proof of Lemma 5.3 of \cite{BD} or
Proposition 3.9 of \cite{B},  shows that
this is impossible.
\end{remark}

\section{Final Remarks}

The arguments given above suggest that there is considerable flexibility in
the construction given by Bourgain and Delbaen. One possibility is to
replace the binary nature of the construction by one which allows a greater
number of terms. Thus in place of $(\pm a,\pm b)$ one might have a
collection of finite sequences $(a_n^j)_{n=1}^N$, $j=1,2,...J.$ Then the
new functionals might evaluate as $\sum_{n=1}^N a_n^k e_{s_n}^*(i_n\pi_n
x-i_{n-1}\pi_{n-1}x)$ where $(s_n)$ is a sequence  such that
$d_{n-1}<s_n\leq d_n$ for each $n$ and $d_n$ is the cardinality of the set
of coordinates defined by the $n$th stage of the construction. Some care
would need to be taken to preserve the boundedness of the iterated
embeddings. It would be most interesting if the set of finite sequences
could be made to vary and if the sequence of finite segments of the
integers could be replaced by finite branches of a tree. This might be an
approach to answering the following question.

\begin{question}
Given a countable ordinal $\alpha$ is there a $\mathcal L_\infty$-space
$X_\alpha$ such that $X_\alpha$ does not contain $c_0$ and $X_\alpha$ has
Szlenk index $\omega^\alpha$?
\end{question}

One other observation is that much of what we have done still works if
$a=1.$ What does not work is the argument in Proposition~\ref{quo}
to find the convergent series for each element
of the dual. Thus the corresponding set $K$ is more complicated and seems
to include a Cantor set of well separated points. A thorough analysis of
this case might yield some additional information about the first example
in \cite{BD}. Finally note that we
have not used the extra conditions imposed on
$a$ and $b$ in \cite{BD} to get a somewhat (hereditarily) reflexive
example.


\begin{thebibliography}{LTII}

\bibitem[A1]{A1}  Dale E. Alspach, \textit{Quotients of }$C[0,1]$
\textit{with
separable dual, }Israel J. Math. 29 (1978), 361--384.
 
\bibitem[A2]{A2}  D. Alspach, \textit{Quotients of }$c_0$\textit{ are
almost
isometric
to subspaces of }$c_0$ , Proc. AMS  76 (1979), 285--288
 
 
\bibitem[A3]{A3}  D. Alspach, \textit{A quotient of $C(\omega^\omega)$
which
is not isomorphic to a subspace of }$C(\alpha),$ $\alpha<\omega_1$,
Israel J. Math  35 (1980), 49--60

\bibitem[A4]{A4}  D. Alspach, \textit{ A }$\ell_1$\textit{-predual
which is not isometric
to a quotient of }$C(\alpha)$, Contemporary Math 144,
Proceedings
of the international workshop in Banach space theory, M\'erida, Venezuela,
(1993), 9--14.

\bibitem[AB]{AB} D. Alspach and Y. Benyamini, $C(K)$\textit{ quotients
of separable }${\mathcal L}_{\infty}$\textit{ spaces},
Israel J. Math 32 (1979), 145--160.

\bibitem[B]{B} J. Bourgain, \textit{New classes
of }${\mathcal {L}}_{p}$\textit{-spaces},
Lecture notes in Mathematics 889, Springer-Verlag, Berlin-New York , 1981.

\bibitem[BD]{BD}  J. Bourgain and F. Delbaen \textit{ A class of special}
$\mathcal L^\infty$\textit{-spaces}, Acta Math  145 (1980), 155--176.

\bibitem[GKL]{GKL} G. Godefroy, N.J. Kalton and G. Lancien, \textit{The
Banach space} $c_0$ \textit{is determined by its metric}, C. R. Acad. Sci.
Paris S\'er. I Math. 327
(1998), 817--822.

\bibitem[GKL1]{GKL1} G. Godefroy, N.J. Kalton and G. Lancien, \textit{Szlenk
indices and uniform homeomorphisms}, preprint.

\bibitem[Ha]{Haydon} R. Haydon, \textit{Subspaces of the Bourgain-Delbaen
space,} preprint.

\bibitem[H]{Haus} F. Hausdorff, \textit{Set Theory}, Second edition. Translated from the German by John R. Aumann et al, Chelsea Publishing Co., New York,
1962.

\bibitem[KOS]{KOS} H. Knaust, E. Odell, and T. Schlumprecht,
\textit{ On asymptotic structure, the Szlenk index and UKK properties in
Banach spaces}, preprint.

\bibitem[JLS]{JLS} W. B. Johnson, J. Lindenstrauss, and G. Schechtman,
\textit{Banach spaces determined by their uniform structures},
 Geom. Funct. Anal. 6 (1996), no. 3, 430--470.

\bibitem[JZ]{JZ} W. B. Johnson and M. Zippin,
{\it Every separable predual of an $L_1$-space is a quotient of }
$C(\Delta)$, Israel J. Math 16 (1973), 198--202.

\bibitem[LS]{LS} D. R. Lewis and C. Stegall, \textit{Banach spaces
whose duals are
isomorphic to }$l_{1}(\Gamma )$, J. Functional Analysis 12 (1973), 177--187.

\bibitem[LTI]{LTI} J. Lindenstrauss and L. Tzafriri, {\it Classical Banach Spaces I, Sequence Spaces}, Springer-Verlag, Berlin, 1977.

\bibitem[LTII]{LTII} J. Lindenstrauss and L. Tzafriri, {\it Classical
Banach Spaces II, Function Spaces,} Springer-Verlag, Berlin, 1979.

\bibitem[S]{S} C. Samuel,
{\it Indice de Szlenk des $C(K)$ ($K$ espace 
topologique compact d\' enombrable)},
Seminar on the  geometry of Banach spaces, Vol. I, II (Paris, 1983),
81--91, Publ. Math. Univ. Paris VII, 18, Univ. Paris VII, Paris, 1984.

\bibitem[Szl]{Szl} W. Szlenk, \textit{The non-existence of a
separable reflexive Banach
space universal for all separable reflexive Banach spaces},
Studia Math.  30 (1968), 53--61.

\bibitem[Z]{Z}  M. Zippin, \textit{The separable extension problem},
Israel J. Math   26 (1977), 372--387.
 
 
\end{thebibliography}
\end{document}